\documentclass[11pt,leqno]{article}

\usepackage{amsmath}
\usepackage{amsfonts}
\usepackage{amssymb}
\date{}
\usepackage{graphicx}
\usepackage{amsmath}

\newcommand{\h}{\hspace*{0.4 cm}}

\newcommand{\ov}{\overline}

\newcommand{\N}{\mathbb{N}}

\newcommand{\R}{\mathbb{R}}
\newcommand{\Z}{\mathbb{Z}}
\newcommand{\E}{{\bf E}}
\newcommand{\F}{{\bf F}}
\newcommand{\G}{{\bf G}}

\title{INTERPOLATION OF BILINEAR OPERATORS AND COMPACTNESS}

\author{\footnotesize E. BRANDANI DA SILVA$^1$ and D. L. FERNANDEZ$^2$}

\begin{document}

\maketitle

{\scriptsize $^1$Universidade Estadual de Maring\'a  - UEM, Departamento de Matem\'atica -
Av. Colombo 5790, Maring\'a - PR, Brazil - 87020-900 \\
\centerline{email: ebsilva@wnet.com.br}}

\vspace{0.7ex}

{\scriptsize $^2$Universidade Estadual de Campinas - Unicamp, Instituto de Matem\'atica -
Caixa Postal 6065, Campinas - SP, Brazil - 13083-859  \\
\centerline{email: dicesar@ime.unicamp.br}}

\vspace{0.7ex}

\h {\footnotesize {\bf Abstract.}
The behavior of bilinear operators acting on interpolation of Banach spaces for the $\rho$ method in relation to the compactness is analyzed. Similar results of Lions-Peetre, Hayakawa and Person's compactness theorems are obtained for the bilinear case and the $\rho$ method.}\\

\vspace{2ex}

{\scriptsize --------------------------------
}\\

\section{Introduction}

The study of compactness of multilinear operators for interpolation spaces goes back to A.P.Calder\'on [4, pp.119-120]. Under an approximation
hypothesis, Calder\'on established an one-side type general result, but
restricted to complex interpolation spaces.

For the real method, if $\E = (E_{0},E_{1}), \F = (F_{0}, F_{1}) \,$ and $\G = (G_{0}, G_{1}) \,$ are Banach couples, a classical result by Lions-Peetre assures that if $T$ is a bounded bilinear operator from $(E_{0} + E_{1}) \times (F_{0} + F_{1}) \,$ into $\, G_{0} + G_{1}, \,$ whose restrictions $T|E_{k} \times F_{k} \; (k =
0,1) \,$ are also bounded from $\, E_{k} \times F_{k} \,$ into $\,
G_{k} \,\, (k = 0,1) \,$, then $T $ is bounded from $\E_{\theta,p;J}
\times \F_{\theta,q;J} \,$ into $\G_{\theta,r;J} \,$, where $\, 0 <
\theta < 1 \,$ and $\, 1/r = 1/p + 1/q - 1 \,$. Lately several authors have obtained new and more general results for interpolation of bilinear and multilinear operators, for example see [9], [12] and [13].

On the other hand, the behaviour of compact multilinear operators under real interpolation functors or more general functors does not seem to have been yet investigated. This is our main subject in this work.

After some preliminaries on interpolation of linear and
bilinear operators, generalizations of Lions-Peetre compactness
theorems [11, Theorem V.2.1] (the one with the same departure spaces) and
[11, Theorem V.2.2] (the one with the same arriving spaces) will be stated. The proof of the first one is an adaptation of the original proof, but the later requires an
involved argument.

Thereafter, a two-side result for general interpolation
functors of type $\rho$, with the additional cost of an approximation
hypothesis on the departure Banach couples, will be then given. Thus, a theorem of Hayakawa type (i.e. a two-side result without approximation hypothesis) will be obtained. The point at this issue is that the Hayakawa type theorem is nothing but a  corollary of the result with approximation hypothesis. Consequently, an one-side result holds for ordered Banach couples.

Finally, as a consequence of the second Lemma of Lions-Peetre type, a compactness theorem of Persson type is obtained. To avoid ponderous notations we have restricted ourselves to
the bilinear case. A generalization for the $\rho$ method of the Lions-Peetre's bilinear theorem will be also provided in this work.

This work was published at "Nonlinear Analysis: Theory, Methods and Applications, Volume 73, Issue 2, 2010, Pages 526-537". Since there are some gaps in the original proof of Theorem 4.3, we give a new proof. For this, we change the Lemma 4.2.

\section{Preliminaries on Interpolation}

{\bf 2.1 Interpolation functors.} A pair of Banach spaces
$\E = (E_{0},E_{1})$ is said to be a Banach couple if $E_{0}$
and $E_{1}$ are continuously embedded in some Hausdorff linear topological space
$\cal E$. Then we can form their intersection $E_{0} \cap E_{1}$ and sum
$E_{0} + E_{1}$; it can be seen that $E_{0} \cap E_{1}$ and
$E_{0} + E_{1}$ become Banach spaces when endowed with the
norms
\[
||\;x\;||_{E_{0} \cap E_{1}} = \max \{\;||x||_{E_{0}},||x||_{E_{1}}\},
\]
and
\[
||\;x\;||_{E_{0}+E_{1}} = \inf_{x=x_{0}+x_{1}} \{||x_{0}||_{E_{0}} +
                                                 ||x_{1}||_{E_{1}}\},
\]
respectively.

\hspace*{4ex}We shall say that a Banach space is an
{\em intermediate space}
in respect to a Banach couple $\E = (E_{0},E_{1})$ if
\[
E_{0} \cap E_{1} \hookrightarrow E \hookrightarrow E_{0}+E_{1}.
\]
(The hookarrow $\hookrightarrow$ denotes bounded embeddings).

\hspace*{4ex}Let $\E = (E_{0},E_{1})$ and $\F = (F_{0},F_{1})$ be
Banach couples. We shall denote by $L(\E,\F)$ the set of all
linear mappings from $E_{0}+E_{1}$ into $F_{0}+F_{1}$ such that
$T|_{E_{K}}$ is {\em bounded} from $E_{k}$ into $F_{k}$, $k=0,1$.

\hspace*{4ex}By an {\em interpolation functor} $\cal F$ we shall mean a functor
which to each Banach couple $\E = (E_{0},E_{1})$ associates an
intermediate space ${\cal F}(E_{0},E_{1})$ between $E_{0}$ and $E_{1}$, and
such that $T|_{{\cal F}(E_{0},E_{1})} \in L({\cal F}(E_{0},E_{1}),
{\cal F}(F_{0},F_{1}))$, for all $T \in L(\E,\F)$.

\hspace*{4ex}The interpolation functors which we shall consider depend
on function parameters.

{\bf 2.2 The function parameters.} By a {\em function parameter} $\rho$ we shall
mean a continuous and positive function on $\R_{+}$.

\hspace*{4ex}We shall say that a function parameter $\rho$ belongs to the
class $\cal B$, if it satisfies the following conditions:
\begin{equation}
   \rho(1) = 1,
\end{equation}
and
\begin{equation}
\overline{\rho}(s) = \sup_{t > 0} \frac{\rho(st)}{\rho(t)} < + \infty,
       \;\;\;\;\; s > 0.
\end{equation}
Also, we shall say that a function parameter $\rho \in \cal B$ belongs to the
class $\cal B^{+-}$ if it satisfies
\begin{equation}
\int_{0}^{\infty} \min (1, \; \frac{1}{t}) \overline{\rho}(t)
       \frac{dt}{t} < + \infty.
\end{equation}
\hspace*{4 ex} From (1)-(3) we see that $\cal B^{+-}$ is contained in
Peetre`s class $\cal P^{+-}$, i.e. the class of pseudo-concave function
parameters which satisfies
\begin{equation}
\overline{\rho}(t) = o (\max (1,t)).
\end{equation}
(See Gustavsson [6] and Gustavsson-Peetre [7].)

\hspace*{4ex}The function parameter $\rho_{\theta}(t) = t^{\theta}$,
$0 \leq \theta \leq 1$, belongs to $\cal B$. It corresponds to the usual
parameter $\theta$. Further, $\rho_{\theta} \in \cal B^{+-}$ if
$0 < \theta < 1$, but $\rho_{0}$, $\rho_{1} \not \in \cal B^{+-}$.

\hspace*{4ex}To control function parameters we shall need to recall the
Boyd indices (see Boyd [2], [3] and Maligranda [12]).

{\bf 2.3 The Boyd indices.} Given a function parameter $\rho \in \cal B$,
the {\em Boyd indices} $\alpha_{\overline{\rho}}$ and
$\beta_{\overline{\rho}}$ of the submultiplicative function $\overline{\rho}$
are defined, respectively, by
\begin{equation}
\alpha_{\overline{\rho}} = \sup_{1 < t < \infty} \frac{\log \overline{\rho}(t)}
                                                      {\log t},
\end{equation}
and
\begin{equation}
\beta_{\overline{\rho}} = \sup_{0 < t < 1} \frac{\log \overline{\rho}(t)}
                                               {\log t}.
\end{equation}
The indices $\alpha_{\overline{\rho}}$ and $\beta_{\overline{\rho}}$ are
real numbers with the followings properties
\begin{equation}
\alpha_{\overline{\rho}} < 0 \;\; \Longleftrightarrow \;\;
  \int_{1}^{\infty} \overline{\rho}(t) \; \frac{dt}{t} < +\infty,
\end{equation}
and
\begin{equation}
\beta_{\overline{\rho}} > 0 \;\; \Longleftrightarrow \;\;
  \int_{0}^{1} \overline{\rho}(t) \; \frac{dt}{t} < +\infty.
\end{equation}

\hspace*{4ex}For the above mentioned function parameters $\rho_{0}$ and
$\rho_{1}$ we have $\alpha_{\overline{\rho}_{0}} < 0$ and
$\beta_{\overline{\rho}_{1}} > 0$, respectively. For the function
parameter $\rho_{\theta}(t) = t^{\theta}$, $0 < \theta < 1$, we have
$\alpha_{\rho_{\theta}} = \beta_{\rho_{\theta}} = \theta.$

\hspace*{4ex}It can be proved that for all $\rho \in \cal B$, with
$\beta_{\overline{\rho}} > 0$ ($\alpha_{\overline{\rho}} < 0$, respectively)
there exists an increasing (decreasing, respectively) function parameter
$\rho^{+}$ ($\rho_{-}$, respectively) equivalent to $\rho$. Hence, if
$\rho \in \cal B^{+-}$ it can be considered an increasing parameter, and
$\rho(t)/t$ a decreasing parameter. Furthermore, $\overline{\rho}$ can be
considered non-decreasing, and $\overline{\rho}(t)/t$ non-increasing.
Consequently, if $\rho \in \cal B^{+-}$ and $0 < q \leq \infty$, we have

\[
||\rho^{-1}(t) \;\min(1,t)\;||_{L^{q}_{*}} < \infty.
\]

{\bf 2.4 Interpolation with function parameters.} Let $\{E_{0},E_{1}\}$ and $\{F_{0},F_{1}\}$
be Banach couples and let $L(\{E_{0},E_{1}\},\{F_{0},F_{1}\})$ be the family
of all linear maps $T : E_{0}+E_{1} \rightarrow F_{0}+F_{1}$ such that
$T|_{E_{k}}$ is {\em bounded} from $E_{k}$ to $F_{k}$, $k=0,1$.

\hspace{4ex}If $E$ and $F$ are intermediate spaces with respect to $\{E_{0},E_{1}\}$
and $\{F_{0},F_{1}\}$, respectively, we say that $E$ and $F$  are interpolation
spaces of type $\rho$ where $\rho \in {\cal P}^{+-}$ if given any $T \in
L(\{E_{0},E_{1}\},\{F_{0},F_{1}\})$ we have

\[
||T||_{L(E,F)} \leq C
||T||_{0} \; \overline{\rho} \biggl( \frac{||T||_{1}}{||T||_{0}} \biggr),
\]
for all $T \in L(\{E_{0},E_{1})\},\{F_{0},F_{1}\})$, where $||T||_{k} = ||T||_{
L(E_{k},F_{k})}$, $(k=0,1)$ and $C > 0$ is a constant.

\hspace{4ex} Let $\{E_{0}, E_{1}\}$ be a Banach couple. The $J$ and $K$
functionals are defined by
\[
J(t,x) = J(t,x;\E) = \max \{||x||_{E_{0}},\; t \;||x||_{E_{1}}\},
                                      \;\;\;\;\ x \in E_{0} \cap E_{1},
\]
\[
K(t,x) = K(t,x;\E) = \inf_{x=x_{0}+x_{1}}
                          \{||x_{0}||_{E_{0}} + \; t \;||x_{1}||_{E_{1}}\},
\]
respectively, where in $K(t,x)$, $x_{0} \in E_{0}$ and $x_{1} \in E_{1}$. Then,
we can define the following interpolation spaces.

\hspace{4ex} The space $(E_{0}, E_{1})_{\rho, q, K}, \;\; \rho \in
{\cal B}$ and $0 < q \leq + \infty$, consists of all $x \in E_{0} + E_{1}$
which norm
\[
|| x ||_{\rho, q; K} = ||(\; \rho(2^{n})^{-1} K(2^{n}, x;\E) \;
)_{n \in \Z}||_{\ell^{q}(\Z)}
\]
is finite.

\hspace{4ex} The space $(E_{0}, E_{1})_{\rho, q; J}$, consists of all $x \in E_{0}
 + E_{1}$, which it has a representation $x = \sum_{n = - \infty}^{\infty} u_{n}$ where $
(u_{n}) \in E_{0} \cap E_{1}$ and converges in $E_{0} + E_{1}$, which
norm
\[
|| x ||_{\rho, q; J} = \inf ||(\; \rho(2^{n})^{-1} J(2^{n},
u_{n}; \E)\;)_{n \in \Z}||_{\ell^{q}(\Z),}
\]
is finite, where the infimum is taken over all representations $x = \sum u_{n}$.
Besides, we have for the interpolation space $(E_{0},E_{1})_{\rho,q,J}$ that
if $x \in E_{0} \cap E_{1}$ then
\[
|| x ||_{E} \leq C || x ||_{0} \;\;
\overline{\rho} \biggl( \frac{||x||_{1}}{||x||_{0}} \biggr).
\]

For $0 < q \leq + \infty$, the Equivalence Theorem between the $J$
and $K$ method holds, that is,
\[
(E_{0},E_{1})_{\rho, q;J} = (E_{0}, E_{1})_{\rho, q;K}.
\]

{\bf 2.5. The spaces of class $J_{\rho}$ and $K_{\rho}$.} Let $E$ be an
intermediate space respect to a Banach couple $\E = (E_{0}, E_{1})$ and $\rho \in {\cal B}^{+-}$. We
say that $E$ is an intermediate space of class $J_{\rho}(E_{0}, E_{1})$ if
the following embedding holds
\begin{equation}
(E_{0}, E_{1})_{\rho, 1; J} \hookrightarrow E,
\end{equation}
and we say that $E$ is an intermediate space of class $K_{\rho}(E_{0},
E_{1})$ if the following embedding holds
\begin{equation}
E \hookrightarrow (E_{0}, E_{1})_{\rho, \infty; K}.
\end{equation}
\hspace*{4ex}We note that $E$ is of class $J_{\theta}(E_{0}, E_{1})$ if and
only if for all $x \in E_{0} \cap E_{1}$, it holds
\begin{equation}
|| x ||_{E} \leq C || x ||_{E_{0}} \ov{\rho}\left(\frac{||x||_{E_{1}}}{|| x ||_{E_{0}}}\right).
\end{equation}

{\bf 2.6 The sequence spaces $\ell_{s}^{q}(G_{m})$}. To obtain our main result, we use the following sequence spaces which are defined as follows.

Let $G$ be a linear space and let $(|| \, \cdot \, ||_{n})_{n \in
\Z}$ be a sequence of norms on $G$.
For each $n \in \Z$, we shall denote by $G_{n}$ the space $G$ equipped
with the norm $||\; \cdot \;||_{n}$:
$G_{n} = (G, || \, \cdot \,||_{n})$.

Let $\rho$ be any function parameter and $0 < q \leq \infty$.
We shall denote by $\ell_{\rho}^{q}(G_{n})$ the linear space of all sequences
$(a_{n})_{n \in  \Z}$, in $G$, such that
\[
|||(a_{n})|||_{\rho, q} = ||(a_{n})_{n \in \Z}||_{\ell_{\rho}^{q}
(G_{n})} = \biggr[ \sum_{n \in \Z} [ \rho(2^{-n})||a_{n}||_{n}]^{q}
\biggr]^{1/q} < + \infty.
\]

The functional $||| \, \cdot \, |||_{\rho ,q}$ is a norm on
$\ell^{q}_{\rho}(G_{n})$. The spaces $\ell_{\rho}^{q}(G_{n})$ are related
with interpolation by the following result:

{\bf Theorem 2.1}. We have for the norm above that
\[
(\ell^{q}_{0}(G_{m}), \ell^{q}_{1}(G_{m}))_{\rho, q} = \ell^{q}_{f}
(G_{m}),          \;\;\;\; 0 < q \leq \infty.
\]
where $f(t) = 1/\rho(t^{-1})$.

For each $m \in \Z$, let us set
\[
\Delta_{m} = \Delta_{m} \E = E_{0} \cap 2^{-m} E_{1} ,
\]
i.e., we take $\Delta_{m}$ to be the space $E_{0} \cap E_{1}$
equipped with the norm $J(2^{-m},\cdot)$.\newline

Giving $\rho \in {\cal B}^{+-}$, for $f(t) = 1/\rho(t)^{-1}$ every sequence $\{u_{m}\}$ in $\ell^{q}_{f}(\Delta_{m})$
is summable in $E_{0} + E_{1}$. Then, setting
\begin{equation}
\sigma(\{u_{m}\}) = \sum^{\infty}_{m=-\infty} u_{m} ,
\end{equation}
by the Theorem 2.1 we see the mapping
\[
\sigma : \ell^{q}_{f} (\Delta_{m}) \longrightarrow (E_{0},
E_{1})_{\rho,q;J}
\]
is bounded and
\[
(E_{0}, E_{1})_{\rho,q;J} = \ell^{q}_{f}(\Delta_{m}) / {\sigma}^{-1} (%
{0}).
\]
Moreover, it can be proved that
\begin{equation}
\ell^{q}_{f}(\Delta_{m}) \subset (\ell^{1}_{0}(\Delta_{m}),
\ell^{1}_{1}(\Delta_{m}))_{\rho,q}. \nonumber
\end{equation}

\section{Bilinear Interpolation}

The following result characterizes the bilinear interpolation
operator which concerns to us. For the classical $\theta$ method
this property was first established
by Lions--Peetre [11, Th.I.4.1]. Here, we give the function parameter version.

Given Banach spaces $E$, $F$ and $G$, we denote by $Bil(E \times
F,G)$ the space of all bilinear operators from $E \times F$ into $G$,
endowed with the norm
\[
||T||_{Bil(X \times Y, Z)} = \sup \{||T(x,y)||_{Z}\; | \;\; || x ||_{X} \leq
1, || y ||_{Y} \leq 1\}.
\]

{\bf Theorem 3.1.} Let $T$ be a bounded bilinear operator from $\, (E_{0} +
E_{1}) \times (F_{0} + F_{1}) \,$ into $\, G_{0} + G_{1} \,$ whose restrictions $\, T|_{E_{k} \times F_{k}} \;\; (k = 0,1) \,$ are bounded from $\, E_{k} \times F_{k} \,$ into $\,
G_{k} \;\; (k = 0,1)$. Then, for $\rho \in {\cal B}^{+-}$ one has
\[
T : E_{\gamma,p} \times F_{\rho,q} \rightarrow G_{\rho,r} \, ,
\]
where $\gamma(t) = \ov{\rho}(t^{-1})^{-1}  \in {\cal B}^{+-}$, $1/r = 1/p + 1/q - 1$ and
\[
||T||_{Bil(E_{\gamma,p} \times F_{\rho,q},G_{\rho,r})} \leq C \, ||T||_{
Bil(E_{0} \times F_{0},G_0)} \, \ov{\rho}\left(\frac{||T||_{
Bil(E_{1} \times F_{1},G_1)}}{||T||_{
Bil(E_{0} \times F_{0},G_0)}}\right) \, ,
\]
where $C > 0$ is a constant.

{\bf Proof.} Let $M_k = ||T||_{
Bil(E_{k} \times F_{k},G_k)}$, $k = 0,1$. Since $\ov{\gamma}(t) \leq \ov{\rho}(t)$ for all $t > 0$, it follows $\gamma(t) \in {\cal B}^{+-}$. Now, let $\, x \in E_{\gamma,p} \,$ and $\, y \in
F_{\rho,q} \,$. Given $\, \varepsilon > 0 \,$, let $\{u_{m}\}$ and $%
\{v_{m}\}$ be sequences in $\, E_{0} \cap E_{1} \,$ and $\, F_{0}
\cap E_{1} \,$, respectively, such that
\[
x = \sum_{m=-\infty}^{\infty} u_{m}\;\; ({\rm in}\; E_{0} +
F_{1}),\;\;\;\; y = \sum_{m=-\infty}^{\infty} v_{m} \;\; ({\rm in}\;
F_{0} + F_{1})
\]
and
\begin{equation}
\left\{
\begin{array}{l}
|| \{\gamma(2^m)^{-1} J(2^{m}, u_{m})\} ||_{\ell^{p}} \leq || x
||_{\E_{\gamma,p}} + \varepsilon  \\
|| \{\rho(2^{m})^{-1} J(2^{m}, v_{m})\} ||_{\ell^{q}} \leq || y
||_{\F_{\rho,q}} + \varepsilon .
\end{array}
\right.
\end{equation}
Hence, for $k = 0,1,$
\begin{equation}
\left\{
\begin{array}{l}
|| \{2^{(k m} \gamma(2^m)^{-1} u_{m} \} ||_{\ell^{p}(E_{k})} \leq || x
||_{\E_{\gamma,p}} + \varepsilon  \\
|| \{2^{(k m} \rho(2^m)^{-1} v_{m} \} ||_{\ell^{q}(F_{k})} \leq || y
||_{\F_{\rho,q}} + \varepsilon .
\end{array}
\right.
\end{equation}
We have
\[
x_{i} = \sum_{| m |\leq i} u_{m} \in E_{0} \cap E_{1} \;\; , \;\;
y_{j} = \sum_{| m |\leq j} v_{m} \in F_{0} \cap F_{1}
\]
and
\[
x_{i} \stackrel{E_{0} + E_{1}}{\longrightarrow} x, \;\; {\rm as}
\;\; i \rightarrow \infty, \;\; y_{j} \stackrel{F_{0} +
F_{1}}{\longrightarrow} y, \;\; {\rm as} \;\; j \rightarrow \infty.
\]
Since,
\begin{eqnarray*}
\lefteqn{|| T(x,y) - T(x_{i}, y_{j}) ||_{G_{0} + G_{1}} \leq} \\
& \leq & || T(x - x_{i},y)||_{G_{0} + G_{1}} + || T(x_{i}, y -
y_{j})
||_{G_{0} + G_{1}} \\
& \leq & M \{|| x - x_{i} ||_{E_{0} + E_{1}} || y ||_{F_{0} + F_{1}}
+ || x_{i} ||_{E_{0} + E_{1}} || y - y_{j} ||_{F_{0} + F_{1}} \} \, ,
\end{eqnarray*}

we see that
\[
T(x, \; y) = \sum_{i \in \Z} \sum_{j \in \Z} T (u_{i},
\; v_{j}) = \sum_{m \in \Z} \sum_{n \in \Z} \;\; T(u_{m}
, \; v_{n - m}) \;\;\; {\rm in} \;\; G_{0} + G_{1}.
\]
Moreover, for each $n \in \Z, \; \sum_{m} T(u_{m}, \; v_{n -
m})$ converges in $G_{0} + G_{1}$, and by Young's inequality and
(13),
\begin{eqnarray}
\lefteqn{||\{2^{k n} \rho(2^n)^{-1} \; \sum_{m} \; ||T(u_{m}, \; v_{n -
m})||_{G_{k}}\}||_{\ell^{r}} } \nonumber \\
& \leq & ||\{2^{k n} \rho(2^n)^{-1} \; \sum_{m} \; M_k ||u_{m}||_{E_k} ||v_{n -
m})||_{F_{k}}\}||_{\ell^{r}} \nonumber \\
& \leq & M_{k} || \{2^{k n} \rho(2^n)^{-1} \; \sum_{m} \; ||u_{m}||_{E_{k}}
\;
||v_{n - m}||_{F_{k}}\}||_{\ell^{r}} \nonumber \\
& = & M_{k} || \{\sum_{m} 2^{k n} 2^{-k m} 2^{k m} \rho(2^n)^{-1} \ov{\rho}(2^{-m})^{-1} \ov{\rho}(2^{-m}) \; ||u_{m}||_{E_{k}} \; ||v_{n - m}||_{F_{k}}\}||_{\ell^{r}} \nonumber \\
& \leq & M_{k} || \{\sum_{m} 2^{k (n - m)} 2^{k m} \rho(2^(n - m))^{-1} \ov{\rho}(2^{-m}) \; ||u_{m}||_{E_{k}} \; ||v_{n - m}||_{F_{k}}\}||_{\ell^{r}} \nonumber \\
& = & M_{k} || \{\sum_{m} 2^{k m} \ov{\rho}(2^{-m}) ||u_{m}||_{E_{k}} 2^{k (n - m)} \rho(2^{(n - m)})^{-1}  \;  ||v_{n - m}||_{F_{k}}\}||_{\ell^{r}} \nonumber \\
& \leq & M_{k}\; || \{ 2^{(k m} \gamma(2^m)^{-1} \;
u_{m}\}||_{\ell^{p}(E_{k})} \; ||\{2^{(k m} \rho(2^m)^{-1} v_{m}\}||_{\ell^{q}(F_{k})} < + \infty \, ,
\end{eqnarray}
for $k = 0, \; 1$. Consequently, for each $n, \; \sum_{m} \; ||T( u_{m}, \; v_{n -
m})||_{G_{0} \cap G_{1}}, \; < + \infty$ and so
\[
w_{m} = \sum_{n} \; T( u_{m}, \; v_{n - m}) \in G_{0} \cap G_{1},
\]
and
\begin{equation}
T(x, \; y) = \sum_{m} w_{m} \;\;\; {\rm in} \;\;\; G_{0} + G_{1} \, .
\end{equation}
From the definition of $\rho$-method for the $J$ functor and from (16) one has,
\begin{eqnarray}
\lefteqn{||T(x, \; y)||_{G_{\rho, \;r}}  \leq \, || \{ \rho(2^{m})^{-1} w_{m}
\}||_{\ell^{r} (G_{0})} \ov{\rho}\left(\frac{||\{ 2^{m} \rho(2^{m})^{-1} \; w_{m} \}
||_{\ell^{r}(G_{1})}}{|| \{ \rho(2^{m})^{-1} w_{m}
\}||_{\ell^{r} (G_{0})}}\right)}  \nonumber \\
& = & C \, \ov{\rho}\left(\frac{1}{\frac{|| \{ \rho(2^{m})^{-1} w_{m}
\}||_{\ell^{r} (G_{0})}}{||\{ 2^{m} \rho(2^{m})^{-1} \; w_{m} \}
||_{\ell^{r}(G_{1})}}} \right)
\left(\frac{1}{||\{ 2^{m} \rho(2^{m})^{-1} \; w_{m} \}
||_{\ell^{r}(G_{1})}} \frac{1}{\frac{|| \{ \rho(2^{m})^{-1} w_{m}
\}||_{\ell^{r} (G_{0})}}{||\{ 2^{m} \rho(2^{m})^{-1} \; w_{m} \}
||_{\ell^{r}(G_{1})}}}\right)^{-1} \, .
\end{eqnarray}

Since $\frac{\ov{\rho}(1/t)}{1/t}$ and $\ov{\rho}(t)$ are non-decreasing, from (15) and (17)

\begin{eqnarray}
\lefteqn{||T(x, \; y)||_{G_{\rho, \;r}}  \leq C \, \ov{\rho}\left(\frac{1}{\frac{M_{0}\; || \{
\gamma(2^m)^{-1} \; u_{m}\}||_{\ell^{p}(E_{0})} \;
||\{\rho(2^m)^{-1} v_{m}\}||_{\ell^{q}(F_{0})}}{||\{ 2^{m}
\rho(2^{m})^{-1} \; w_{m} \} ||_{\ell^{r}(G_{1})}}} \right)} \nonumber \\
& & \times \left(\frac{1}{||\{ 2^{m} \rho(2^{m})^{-1} \; w_{m} \}
||_{\ell^{r}(G_{1})}} \frac{1}{\frac{M_{0}\; || \{ \gamma(2^m)^{-1}
\; u_{m}\}||_{\ell^{p}(E_{0})} \; ||\{\rho(2^m)^{-1}
v_{m}\}||_{\ell^{q}(F_{0})}}{||\{ 2^{m} \rho(2^{m})^{-1} \; w_{m} \}
||_{\ell^{r}(G_{1})}}}\right)^{-1} \nonumber \\
& \leq & C \, M_0 \, \ov{\rho}\left(\frac{||\{
2^{m}\rho(2^{m})^{-1} \; w_{m} \} ||_{\ell^{r}(G_{1})}}{M_0 \, || \{
\gamma(2^m)^{-1} \; u_{m}\}||_{\ell^{p}(E_{0})} \;
||\{\rho(2^m)^{-1} v_{m}\}||_{\ell^{q}(F_{0})}} \right) \nonumber \\
& & \times \left(\frac{1}{|| \{ \gamma(2^m)^{-1} \; u_{m}\}||_{\ell^{p}(E_{0})} \;
||\{\rho(2^m)^{-1} v_{m}\}||_{\ell^{q}(F_{0})}}\right)^{-1} \nonumber \\
& \leq & C \, M_0 \, \ov{\rho}\left(\frac{M_{1}\; || \{ 2^{m}
\gamma(2^m)^{-1} \; u_{m}\}||_{\ell^{p}(E_{1})} \; ||\{2^{m}
\rho(2^m)^{-1} v_{m}\}||_{\ell^{q}(F_{1})}}{ M_{0}\; || \{
\gamma(2^m)^{-1} \; u_{m}\}||_{\ell^{p}(E_{0})} \;
||\{\rho(2^m)^{-1} v_{m}\}||_{\ell^{q}(F_{0})}} \right) \nonumber \\
& & \times \left(\frac{1}{
|| \{ \gamma(2^m)^{-1} \; u_{m}\}||_{\ell^{p}(E_{0})} \;
||\{\rho(2^m)^{-1} v_{m}\}||_{\ell^{q}(F_{0})}}\right)^{-1} \, .
\end{eqnarray}

Finally, since $\ov{\rho}(1/t)(1/t)^{-1}$ and $\ov{\rho}(t)$ are submultiplicative, from (18) one has
\begin{eqnarray*}
\lefteqn{||T(x, \; y)||_{G_{\rho, \;r}} \leq C \, M_0 \, \ov{\rho}(\frac{M_1}{M_0})
\ov{\rho}\left(\frac{1}{ || \{ \gamma(2^m)^{-1} \;
u_{m}\}||_{\ell^{p}(E_{0})} \; ||\{\rho(2^m)^{-1}
v_{m}\}||_{\ell^{q}(F_{0})}} \right)} \\
& & \times \, \left(\frac{1}{ || \{
\gamma(2^m)^{-1} \;
u_{m}\}||_{\ell^{p}(E_{0})} \; ||\{\rho(2^m)^{-1} v_{m}\}||_{\ell^{q}(F_{0})}}\right)^{-1}
\\
& & \times \, \ov{\rho}(
|| \{ 2^{m} \gamma(2^m)^{-1} \;
u_{m}\}||_{\ell^{p}(E_{1})} \; ||\{2^{m} \rho(2^m)^{-1} v_{m}\}||_{\ell^{q}(F_{1})})  \\
& \leq & C \, M_0 \, \ov{\rho}(\frac{M_1}{M_0})
\ov{\rho}\left(\frac{1}{ || \{ \gamma(2^m)^{-1} \;
u_{m}\}||_{\ell^{p}(E_{0})}} \right) \left(\frac{1}{ || \{
\gamma(2^m)^{-1} \; u_{m}\}||_{\ell^{p}(E_{0})}}\right)^{-1} \\
& & \times \, \ov{\rho}\left(\frac{1}{ ||\{\rho(2^m)^{-1}
v_{m}\}||_{\ell^{q}(F_{0})}}\right) \left(\frac{1}{
||\{\rho(2^m)^{-1} v_{m}\}||_{\ell^{q}(F_{0})}}\right)^{-1} \\
& & \times \, \ov{\rho}(
|| \{ 2^{m} \gamma(2^m)^{-1} \;
u_{m}\}||_{\ell^{p}(E_{1})}) \ov{\rho}(||\{2^{m} \rho(2^m)^{-1} v_{m}\}||_{\ell^{q}(F_{1})})  \\
& \leq & C \, M_0 \, \ov{\rho}(\frac{M_1}{M_0})
\ov{\rho}\left(\frac{1}{|| x ||_{\E_{\gamma,p}} + \varepsilon}
\right) \left(\frac{1}{|| x ||_{\E_{\gamma,p}} + \varepsilon}\right)^{-1}
\ov{\rho}\left(\frac{1}{|| y ||_{\F_{\rho,q}} +
\varepsilon} \right)\left(\frac{1}{|| y ||_{\F_{\rho,q}} +
\varepsilon}\right)^{-1} \\
& & \times \, \ov{\rho}(|| x ||_{\E_{\gamma,p}} + \varepsilon)
\ov{\rho}(|| y ||_{\F_{\rho,q}} + \varepsilon) \, .
\end{eqnarray*}
Since $\varepsilon$ is arbitrary, for $x$ and $y$ fixed, taking $\varepsilon \rightarrow 0$ one has
\begin{eqnarray*}
||T(x, \; y)||_{G_{\rho, \;r}} &\leq& M_0 \ov{\rho}(\frac{M_1}{M_0})
\ov{\rho}\left(\frac{1}{|| x ||_{\E_{\gamma,p}}}
\right) \left(\frac{1}{|| x ||_{\E_{\gamma,p}}}\right)^{-1}
\ov{\rho}\left(\frac{1}{|| y ||_{\F_{\rho,q}}}
\right)\left(\frac{1}{|| y ||_{\F_{\rho,q}}}\right)^{-1} \\
 & & \times \, \ov{\rho}(|| x ||_{\E_{\gamma,p}})
\ov{\rho}(|| y ||_{\F_{\rho,q}}) \, .
\end{eqnarray*}

Since $\ov{\rho}(1) = 1$, taking the supremum on $|| x ||_{\E_{\gamma,p}} \leq 1$ and $|| y ||_{\F_{\rho,q}} \leq 1$ the result follows.

For more on bilinear interpolation see Janson [5], Maligranda [8] and Mastylo [9].

\section{Compactness Theorems of Lions-Peetre Type}

Given Banach spaces $E, F$ and $G$, a bounded bilinear mapping $T$ from $E \times F$ into $%
G$ is {\em compact} if the image of the set $M = \{(x,y) \in E \times F : \max\{||x||_{E}, \; ||y||_{F}\} \leq 1 \}$ is a totally bounded subset of $G$.

In this section versions of Lions-Peetre Theorems about compactness of operators in interpolated spaces by the $\rho$ method shall be given. We begin with the version of the Theorem V.2.1([11]).

{\bf Theorem 4.1.} Let $E$ and $F$ be Banach spaces, $\G =
(G_{0}, G_{1})$ a Banach couple and $G$ be a Banach space of class $%
J_{\rho} (G_{0}, G_{1})$, $\rho \in {\cal B}^{+-}$. Given a bounded
bilinear operator $T$ from $E \times F$ into $G_{0} + G_{1}$, such that $T(E
\times F) \subset G_{0} \cap G_{1}$, $T$ is compact
from $E \times F$ into $G_{0}$ and bounded from $E \times F$ into $G_{1}$, then $T$ is also compact from $E \times F$ into $G$.

{\bf Proof.} Let $\{(x_{n}, y_{n}) \}$ be a bounded sequence in $E \times
F$. Since $T$ is compact from $E \times F$ into $G_{0}$, there exists a
subsequence $\{(x_{\nu} , y_{\nu}) \}$ such that $\{T( x_{\nu} , y_{\nu}) \}$
is a $G_{0}$-Cauchy sequence. On the other hand, since $G$ belongs to class $%
J_{\rho} (G_{0} , \; G_{1})$, we have

\begin{eqnarray*}
\lefteqn{||T(x_{\mu} , y_{\mu}) - T(x_{\nu} , y_{\nu})||_{G} \leq} \\
& \leq & C \, ||T(x_{\mu} , y_{\mu}) - T(x_{\nu} , y_{\nu})||_{G_0} \; \ov{\rho}\left(\frac{||T(x_{\mu} , y_{\mu}) - T(x_{\nu} , y_{\nu})
||_{G_{1}}}{||T(x_{\mu} , y_{\mu}) - T(x_{\nu} , y_{\nu})||_{G_0}}\right) \\
& \leq & C \, \ov{\rho}(||T||_{Bil(E \times F, \; G_{1})}) \, \frac{\ov{\rho}\left(\frac{1}{||T(x_{\mu} , y_{\mu}) - T(x_{\nu} , y_{\nu})||_{G_0}}\right)}{\frac{1}{||T(x_{\mu} , y_{\mu}) - T(x_{\nu} , y_{\nu})||_{G_0}}}.
\end{eqnarray*}

Hence $\{T(x_{\nu} , y_{\nu})\}$ is a $G$-Cauchy sequence and $T$ is compact
from $E \times F$ into $G$.

Now, we turn over the ``Lions-Peetre Theorem V.2.2 ([11])".
A preliminary result which depends on an approximation
hypothesis is previously required.

{\bf Approximation Hypothesis}. A Banach couple ${\cal X} = (X_{0}, X_{1})$ verifies the approximation hypothesis (AP) if there exists a
sequence $\{P_{n}\}$ in $L({\cal {X}},{\cal {X}})$, with $P_{n}(X_{0} + X_{1}) \subset
X_{0} \cap X_{1}$, and two other sequences $\{P^{+}_{n}\}$ and $\{P^{-}_{n}\}
$ in $L({\cal {X}},{\cal {X}})$, such that

{\bf (AP1)} They are uniformly bounded in $L({\cal {X})}$;

{\bf (AP2)} $I = P_{n} + P^{+}_{n} + P^{-}_{n}, \;\;\; n \in \N$ ;

{\bf (AP3)} $P^{+}_{n} = P^{+}_{n} |_{X_{0}} \in L(X_{0}, \; X_{1})$ and $%
P^{-}_{n} = P^{-}_{n} |_{X_{1}} \in L(X_{1}, \; X_{0})$,
and
\[
\lim_{n \rightarrow \infty} ||P^{+}_{n}||_{L(X_{0}, \; X_{1})} = \lim_{n
\rightarrow \infty} ||P^{-}_{n}||_{L(X_{1}, \; X_{0})} = 0.
\]

The following lemma is required.

{\bf Lemma 4.2.} Given Banach spaces $E, F$ and $G$, let $T_n \in Bil(E \times F,G)$, $n \in \N$ be a sequence of bilinear operators such that $\lim_{n \rightarrow \infty} \|T_n\|_{Bil(E \times F,G)} = \lambda$. Then, there exits a sequence $(x_n,y_n) \subset U_E \times U_F$ such that
\[
\lim_{n \rightarrow \infty} \|T_n(x_n,y_n)\|_{G} = \lambda \, .
\]

{\bf Proof}. Since $\lim_{n \rightarrow \infty} \|T_n\|_{Bil(E \times F,G)} = \lambda$, given $\varepsilon > 0$ there is $N > 0$ with $|\|T_n\|_{Bil(E \times F,G)} - \lambda| < \varepsilon$ for all $n > N$. Then,
\[
\|T_n(x,y)\|_G \leq \|T_n\|_{Bil(E \times F,G)} < \varepsilon + \lambda
\]
for all $(x,y) \in U_E \times U_F$. On other hand, since $\|T_n\|_{Bil(E \times F,G)} = \sup_{(x,y) \in U_E \times U_F} \|T_n(x,y)\|_G$, by the definition of supremum, for the given $\varepsilon > 0$ there is $(x_n,y_n) \in U_E \times U_F$, such that $\|T_n\|_{Bil(E \times F,G)} - \varepsilon < \|T_n(x_n,y_n)\|_G \leq \|T_n\|_{Bil(E \times F,G)} < \varepsilon + \lambda$, this implies $\|T_n\|_{Bil(E \times F,G)} - \varepsilon -\lambda < \|T_n(x_n,y_n)\|_G - \lambda < \varepsilon$, since $-\varepsilon < \|T_n\|_{Bil(E \times F,G)} - \lambda$. Thus, one has $-2 \varepsilon < \|T_n(x_n,y_n)\|_G - \lambda < \varepsilon < 2 \varepsilon$, for all $n > N$, which gives the desired limit.

{\bf Theorem 4.3.} Let $\E = (E_0, F_1)$ and $\F = (F_0, F_1)$ be
Banach couples satisfying the Approximation Hypothesis. Given $\rho \in {\cal B}^{+-}$, a Banach space $G$ and $T \in Bil((E_0 + E_1) \times (F_0 + F_1), G)$, such that $T$ is compact from $E_0 \times F_0$ into $G$,
then $T$ is also compact from $\E_{\gamma, p} \times \F_{\rho, q}$ into $G$, for $1 \leq p, q \leq \infty$ and $\gamma(t) = 1/\ov{\rho}(t^{-1})$.

{\bf Proof. Step 1:} Let $(P_m)$ and $(Q_m)$ be the approximating
sequences in $\E$ and $\F$, respectively. We are considering here
the following scheme:

\[
\E_{\gamma, p} \times \F_{\rho, q} \; {\stackrel{\scriptstyle{%
(P_m, Q_m)}}{\longrightarrow}} E_0 \cap E_1 \times F_0 \cap F_1 {\stackrel{\scriptstyle{%
i_j}}{\longrightarrow}} E_j \times F_j {\stackrel{\scriptstyle{%
T}}{\longrightarrow}} G,
\]

for $j = 0,1$ where $i_j$ is the inclusion operator. Then, $T (P_m, Q_m)$ is compact from $%
\E_{\gamma, q} \times \F_{\rho, q}$ into $G$. To prove the
compactness of $T$ it must be shown that
\[
|| T - T(P_m, Q_m) ||_{Bil (\E_{\gamma, p} \times \F_{\rho, q}, G)}
\longrightarrow 0 \;\;\; {\rm as} \;\;\; m \longrightarrow \infty.
\]
Since,
\begin{eqnarray*}
T - T(P_m, Q_m) & = & T(I, I) - T(I, Q_m) - T(P_m, Q_m) + T(P_m,
Q_m) + \\ \nonumber & \;\; & T(I, Q_m) - T(P_m, Q_m) + T(P_m, I) -
T(P_m, Q_m) \\ \nonumber & = & T(I - P_m, I - Q_m) + T(I - P_m, Q_m)
+ T(P_m, I - Q_m) \\ \nonumber & = & T(P^{+}_m, + P^{-}_m, Q^{+}_m +
Q^{-}_m) + T(P^{+}_m, + P^{-}_m, Q_m) + T(P_m, Q^{+}_m + Q^{-}_m) \\
\nonumber & = & T(P^{+}_m, Q^{+}_m) + T(P^{+}_m, Q^{-}_m) +
T(P^{-}_m, Q^{+}_m) + T (P^{-}_m, Q^{-}_m) + \\ \nonumber & \;\; &
T(P^{+}_m, Q_m) + T(P^{-}_m, Q_m) + T(P_m, Q^{+}_m) + T(P_m,Q^{-}_m)
\,\,\,\,\, (*),
\end{eqnarray*}

it needs be shown that each term in $(*)$ goes to zero in
\mbox{$Bil(\E_{\gamma, p} \times \F_{\rho, q}, G)$}.

{\bf Step 2:} We will show that $\displaystyle{\lim_{m \rightarrow
\infty}} ||T(P^{+}_m, Q^{+}_m) ||_{Bil (\E_{\gamma, q} \times
\F_{\rho, q}, G)} = 0$.

By the hypothesis we have
\[
||T(P^{+}_{n},Q^{+}_{n})||_{Bil(\E_{\gamma,p} \times \F_{\rho,q},G)} \leq C \, ||T(P^{+}_{n},Q^{+}_{n})||_{0} \, \ov{\rho}\left(\frac{||T(P^{+}_{n},Q^{+}_{n})||_{1}}{||T(P^{+}_{n},Q^{+}_{n})||_{0}}\right),
\]

where $|| . ||_{k} = ||.||_{Bil(E_{k} \times F_{k},G)}$, $k = 0,1$. Let $x \in B_{E_0}$ and $y \in B_{F_0}$, since $P_{n}^{+} : E_{0} \rightarrow E_{0} \cap E_{1}$ and $Q_{n}^{+} : F_{0} \rightarrow F_{1}$, one has
\[
E_{0} \times F_{0} \stackrel{(P_m^{+}, Q_m^{+})}{\longrightarrow} (E_{0} \cap
E_{1}) \times F_{1} \stackrel{T}{\longrightarrow} G
\]
and
\begin{eqnarray*}
||T||_{Bil((E_{0} \cap E_{1}) \times F_{k})} &=& \sup_{(a,b) \in B_{E_{0} \cap E_{1}} \times B_{F_{k}}} ||T(a,b)||_{G} \\
& \leq & \sup_{(a,b) \in B_{E_k} \times B_{F_k}} ||T(a,b)||_{G)} \\
& = & ||T||_{Bil(E_k \times F_k,G)} \,\,\,\,\, , \,\,\, k = 0,1 \, ,
\end{eqnarray*}

where $B_E$ denotes the unitary ball in the Banach space $E$. Thus,
\begin{eqnarray*}
||T(P^{+}_m x, Q^{+}_m y)||_{G} & \leq & C ||P^{+}_m x||_{E_0 \cap E_1}
||Q^{+}_m y||_{F_1} \\
& \leq & C (||P^{+}_m x ||_{E_0} + ||P^{+}_m x ||_{E_1}) || Q^{+}_m y||_{F_1}
\\
& \leq & C (||P^{+}_m ||_{L (E_0, E_0)} + ||P^{+}_m||_{L(E_0, E_1)})
||Q^{+}_m ||_{L (F_0, F_1)} || x ||_{E_0} ||y||_{F_0} \\
& \leq & C ||Q^{+}_m ||_{L (F_0, F_1)}.
\end{eqnarray*}

By $(AP3)$, $||Q^{+}_m ||_{L (F_0, F_1)} \rightarrow 0$ when $n \rightarrow \infty$, then we get
$||T(P^{+}_m x, Q^{+}_m y)||_{0} \rightarrow 0$ when $n \rightarrow \infty$, which proves the assertion.

{\bf Step 3:} \hspace{4ex} $\displaystyle{\lim_{m \rightarrow
\infty}} || T(P^{+}_m, Q^{-}_m)||_{Bil (\E_{\gamma, q} \times
\F_{\rho q}, G)} = 0$.

We have
\[
|| T(P^{+}_m, Q^{-}_m)||_\theta \leq C \, || T(P^{+}_m, Q^{-}_m)||_0 \, \ov{\rho}\left(\frac{|| T(P^{+}_m, Q^{-}_m)||_1}{||T(P^{+}_{n},Q^{-}_{n})||_{0}}\right) .
\]
Both factors to the right in the above inequality are bounded. Checking that
\[
\lim_{m \rightarrow \infty} || T(P^{+}_m, Q^{-}_m)||_{Bil(E_{0} \times F_{0},G)} = 0 \, ,
\]
is sufficient.

Suppose $|| T(P^{+}_m, Q^{-}_m)||_{Bil(E_{0} \times F_{0},G)} \nrightarrow 0$. Then, there exists $\lambda > 0$ and a subsequence $\{m'\}$ such that $|| T(P^{+}_{m'}, Q^{-}_{m'})||_{Bil(E_{0} \times F_{0},G)} > \lambda$ for all $\{m'\}$.

Since the sequence $\{T(P^{+}_m, Q^{-}_m)\}$ is uniformly bounded, and
\[
0 < \lambda < || T(P^{+}_{m'}, Q^{-}_{m'})||_{Bil(E_{0} \times F_{0},G)} \leq ||T||_{Bil(E_{0} \times F_{0},G)} \|P^+_m\|_{L(E_0,E_0)} \|Q^-_m\|_{L(E_0,E_0)}
\]
the sequence $\{|| T(P^{+}_{m'}, Q^{-}_{m'})||_{Bil(E_{0} \times F_{0},G)}\}_{m'}$ has a convergent subsequence

$\{|| T(P^{+}_{m''}, Q^{-}_{m''})||_{Bil(E_{0} \times F_{0},G)}\}_{m''}$ to $k \geq \lambda > 0$. By the Lemma 4.2 there is a sequence $(x_{m''},y_{m''}) \in U_{E_0} \times U_{F_0}$ such that $\lim_{m'' \rightarrow \infty} || T(P^{+}_{m''}, Q^{-}_{m''})||_{Bil(E_{0} \times F_{0},G)} = k > 0$. By the hypothesis, $T \in Bil((E_0 + E_1) \times (F_0 + F_1),G)$ and since $P^+_n : E_0 \rightarrow E_0 + E_1$ and $Q_n^- : F_0 \rightarrow F_0 + F_1$, considering that
\[
E_{0} \times F_{0} \stackrel{(P_{m}^{+},Q_{m}^{-})}{\longrightarrow} (E_{0} +
E_{1}) \times (F_0 + F_{1}) \stackrel{T}{\longrightarrow} G,
\]
one has
\begin{eqnarray*}
||T(P^{+}_{m^{\prime \prime}} x_{m''}, Q^{-}_{m^{\prime \prime}} y_{m''})||_G & \leq & ||T||_{Bil((E_0 + E_1) \times (F_0 + F_1),G)} ||P^{+}_{m^{\prime \prime}} x_{m''}||_{E_{0} + E_1} || Q^{-}_{m^{\prime \prime}} y_{m''}||_{F_0 + F_1} \\
& \leq & C \; ||P^{+}_{m^{\prime \prime}} x_{m''} ||_{E_1} ||Q^{-}_{m^{\prime \prime}} y_{m''}||_{F_0} \\
& \leq & C \; ||P^{+}_{m^{\prime \prime}}||_{L(E_0,E_1)} || x_{m''}||_{E_{0}} ||Q^{-}_{m^{\prime \prime}}||_{L(F_0,F_0)} ||y_{m''} ||_{F_0} \, ,
\end{eqnarray*}
and, by the hypothesis (AP), $||P^{+}_{m^{\prime \prime}}||_{L(E_0,E_1)} \rightarrow 0$ for $m^{\prime \prime} \rightarrow \infty$. Thus, we have a contradiction.

{\bf Step 4:} $\displaystyle{\lim_{m \rightarrow \infty}}
||T(P^{-}_m, Q^{-}_m)||_{Bil (\E_{\gamma, q} \times \F_{\rho, q},
G)} = 0$.

The proof is symmetrical to step 2, considering the scheme
\[
E_{1} \times F_{1} \stackrel{(P_m^{-}, Q_m^{-})}{\longrightarrow} (E_{0} \cap
E_{1}) \times F_{0} \stackrel{T}{\longrightarrow} G
\]
{\bf Step 5:} $\displaystyle{\lim_{m \rightarrow \infty}} || T(P^{-}_m,
Q^{+}_m)||_{Bil (\E_{\gamma, q} \times \F_{\rho q}, G)} = 0$

The proof is symmetrical to step 3, considering the scheme
\[
E_{0} \times F_{0} \stackrel{(P_{m}^{-},Q_{m}^{+})}{\longrightarrow} (E_{0} +
E_{1}) \times (F_0 + F_{1}) \stackrel{T}{\longrightarrow} G \, .
\]
{\bf Step 6:} $\displaystyle{\lim_{m \rightarrow 0}} ||
T(P_m, Q^{+}_m)||_{Bil (\E_{\gamma, q} \times \F_{\rho, q}, G)} =
\displaystyle{\lim_{m \rightarrow \infty}} ||T(P_m, Q^{-}_m)||_{Bil
(\E_{\gamma, q} \times \F_{\rho, q}, G)} = 0$.

For $T(P_m,Q^{+}_m)$ one has
\[
|| T(P_m, Q^{+}_m)||_{Bil (\E_{\gamma q} \times \F_{\rho, q}, G)}
\leq C \, || T(P_m, Q^{+}_m)||_0 \, \ov{\rho}\left(\frac{|| T(P_m,
Q^{+}_m)||_1}{|| T(P_m, Q^{+}_m)||_0}\right) .
\]
Both factors to the right in the above inequality are bounded. Checking that the former tends to $0$ is sufficient.

Suppose $|| T(P_m, Q^{+}_m)||_{Bil(E_{0} \times F_{0},G)} \nrightarrow 0$. Then, there exists $\lambda_0 > 0$ and a subsequence $\{m'\}$ such that $|| T(P_{m'}, Q^{+}_{m'})||_{Bil(E_{0} \times F_{0},G)} > \lambda_0$ for all $\{m'\}$.

Since the sequence $\{T(P_m, Q^{+}_m)\}$ is uniformly bounded, and
\[
0 < \lambda_0 < || T(P_{m'}, Q^{+}_{m'})||_{Bil(E_{0} \times F_{0},G)} \leq ||T||_{Bil(E_{0} \times F_{0},G)} \|P_m\|_{L(E_0,E_0)} \|Q^+_m\|_{L(E_0,E_0)}
\]
the sequence $\{|| T(P_{m'}, Q^{+}_{m'})||_{Bil(E_{0} \times F_{0},G)}\}_{m'}$ has a convergent subsequence

$\{|| T(P_{m''}, Q^{+}_{m''})||_{Bil(E_{0} \times F_{0},G)}\}_{m''}$ to $k \geq \lambda > 0$. By the Lemma 4.2 there is a sequence $(x_{m''},y_{m''}) \in U_{E_0} \times U_{F_0}$ such that $\lim_{m'' \rightarrow \infty} || T(P_{m''}, Q^{+}_{m''})||_{Bil(E_{0} \times F_{0},G)} = k > 0$. By the hypothesis, $T \in Bil((E_0 + E_1) \times (F_0 + F_1),G)$, then
\begin{eqnarray*}
\|T\|_{Bil(E_0 \times (F_0 + F_1),G)} & = & \sup_{(x,y) \in U_{E_0} \times U_{F_0 + F_1}} \|T(x,y)\|_G \\
& \leq & \sup_{(x,y) \in U_{E_0 + E_1} \times U_{F_0 + F_1}} \|T(x,y)\|_G \\
& = & \|T\|_{Bil((E_0 + E_1) \times (F_0 + F_1),G)} \, ,
\end{eqnarray*}
and since $P_n : E_0 \rightarrow E_0$ and $Q_n^+ : F_0 \rightarrow F_0 + F_1$, considering that
\[
E_{0} \times F_{0} \stackrel{(P_{m},Q_{m}^{+})}{\longrightarrow} E_{0} \times (F_0 + F_{1}) \stackrel{T}{\longrightarrow} G,
\]
one has
\begin{eqnarray*}
||T(P_{m^{\prime \prime}} x_{m''}, Q^{+}_{m^{\prime \prime}} y_{m''})||_G & \leq & ||T||_{Bil(E_0 \times (F_0 + F_1),G)} ||P_{m^{\prime \prime}} x_{m''}||_{E_{0}} || Q^{+}_{m^{\prime \prime}} y_{m''}||_{F_0 + F_1} \\
& \leq & C \; ||P_{m^{\prime \prime}} x_{m''} ||_{E_0} ||Q^{+}_{m^{\prime \prime}} y_{m''} ||_{F_1} \\
& \leq & C \; ||P_{m^{\prime \prime}}||_{L(E_0,E_0)} || x_{m''}||_{E_{0}} ||Q^{+}_{m^{\prime \prime}}||_{L(F_0,F_1)} ||y_{m''} ||_{F_0} \, ,
\end{eqnarray*}
and $||Q^{+}_{m^{\prime \prime}}||_{L(F_0,F_1)} \rightarrow 0$ for $m^{\prime \prime} \rightarrow \infty$, by the hypothesis AP. Then we get a contradiction.

Using a similar reasoning, considering the scheme
\[
E_{1} \times F_{1} \stackrel{(P_{m},Q_{m}^{-})}{\longrightarrow} E_{1} \times (F_0 + F_{1}) \stackrel{T}{\longrightarrow} G,
\]
it is obtained that
$\displaystyle{\lim_{m \rightarrow \infty}}
||T(P_m, Q^{-}_m)||_{Bil(E_1 \times F_1,G)} = 0$.

The proof is complete.

The second theorem of
Lions--Peetre type may be established now.

{\bf Theorem 4.3.} Let $E$ and $F$ be Banach spaces
of class $K_{\gamma}(E_0, E_1)$ and $K_{\rho}(F_0, F_1)$ respectively, where $\rho \in {\cal B}^{+-}$ and $\gamma(t) = 1/\ov{\rho}(t^{-1})$ and let $G$ be any
Banach space. Given $T \in Bil((E_0 + E_1) \times (F_0 + F_1), G)$, such that $T$ is compact from $E_0 \times F_0$ into $G$, then $T$ is also compact from $E \times F$ into $G$. \newline

{\bf Proof.} Let $\sigma$ be the mapping given by (12). We shall prove
that
\begin{equation}
\widetilde{T} = T \; \circ \; (\sigma, \sigma) : \ell^{\infty}_\gamma
(\Delta_m \E) \times \ell^{\infty}_\rho (\Delta_m \F) \longrightarrow
G
\end{equation}

is compact. For either $X = \E$ or $X = \F$, we have
\[
\overline{\ell^{1}_0 (\Delta_m X) \cap \ell^{1}_1 (\Delta_m X)}^k =
\ell^{1}_k (\Delta_m X), \;\; k = 0, 1.
\]

For each $n \in \Z$, let us consider the cutting operators $\, P_{n}
, P^{+}_{n} \,$ and $\, P^{-}_{n} \,$, defined on $\,
\ell^{1}_{0}(\Delta_{m}X) + \ell^{1}_{1}(\Delta_{m}X) \,$ by
\begin{eqnarray*}
P_{n} ({u_{m}}) & = & \{\cdots, 0, 0, u_{-n}, u_{-n+1} \cdots, u_{0},\cdots,
u_{n-1}, u_{n}, 0, 0,\cdots \}, \\
P^{+}_{n} ({u_{m}}) & = & \{ \cdots, 0, 0, u_{n+1}, u_{n+2},\cdots \}, \\
P^{-}_{n} ({u_{m}}) & = & \{ \cdots, u_{-n-2}, u_{-n-1}, 0, 0,\cdots \},
\end{eqnarray*}
thus, the Banach couple $(\ell^{1}_0 (\Delta_m X), \ell^{1}_1 (\Delta_m X))$
verifies the Approximation Hypothesis, and
\[
\ell^{\infty}_f (\Delta_m X) \subset (\ell^{1}_0 (\Delta_m X), \ell^{1}_1
(\Delta_m X))_{\rho, \infty} ,
\]
where either $f(t) = 1/\rho(t)^{-1}$ or $f(t) = 1/\gamma(t)^{-1}$. Since the conditions of Theorem 4.3 are verified, the
bilinear mapping (20) is compact and, {\em a fortiori}, $T$ is also compact.

\section{Compactness Theorem of Hayakawa Type}

In this section we establish a bilinear version of Hayakawa's compactness theorem [8] for the $\rho$ method, in which we assume compactness in both departure spaces and any inclusion conditions. We begin with an preliminary result.

{\bf Theorem 5.1.} Let us assume that $\E = (E_{0}, E_{1})$ and $\F
= (F_{0}, F_{1})$ are Banach couples which satisfy the approximation
hypothesis (AP) and that $\G = (G_{0}, G_{1})$ is an arbitrary
Banach couple. Let $T \in Bil((E_0 + E_1) \times (F_0 + F_1), G_0 + G_1)$ such that the
restrictions $T|_{E_{k} \times F_{k}} \; (k = 0, 1)$ are bounded and
compact from $E_{k} \times F_{k}$ into $G_{k}$, $(k = 0, 1)$. Given $\rho \in {\cal B}^{+-}$ and $\gamma(t) =1/\ov{\rho}(t^{-1})$, if $E \in J_\gamma(E_0,E_1) \cap K_\gamma(E_0,E_1)$, $F \in J_\rho(F_0,F_1) \cap K_\rho(F_0,F_1)$ and $G \in J_\rho(G_0,G_1) \cap K_\rho(G_0,G_1)$, then $T$ is
also compact from $E \times F$ into $G$. Moreover, if $E_{0}
\hookrightarrow E_{1}$ and $F_{0} \hookrightarrow F_{1}$, it is
enough to consider compactness only from $E_{0} \times F_{0}$ into
$G_{0}$.

{\bf Proof. Step 1:} Let $\{P_{m}\}, \; \{P^{+}_{m}\}, \;
\{P^{-}_{m}\}$ and $\{Q_{n}\}, \; \{Q^{+}_{n}\}, \; \{Q^{-}_{n}\}$
be approximating sequences in the Banach couples $\E$ and $\F$,
respectively, satisfying the Approximation Conditions (AP1)-(AP3).

To show that $T : E \times F \rightarrow G$ is compact, it is
suffices to prove that:

i) $T \circ (P_{m}, Q_{n}), \; T \circ (P_{m}, I)$ and $T \circ (I, Q_{n})$
are compact from $E \times F$ into $G$.

ii) $\lim_{n \rightarrow \infty} ||T - T \circ (P_{m}, I) - T \circ
(I,Q_{n}) + T \circ (P_{m},Q_{n})||_{Bil(E \times F,G)} = 0$.

To prove i) we factorize $T \circ (P_{m},Q_{n})$ using the
following diagram:
\[
E \times F \stackrel{(P_{m},Q_{n})}{\longrightarrow} (E_{0} \cap E_{1})
\times (F_{0} \cap F_{1}) \hookrightarrow E_{j} \times F_{j} \stackrel{T}{\longrightarrow} G_{j},
\]
for $j = 0,1$.

Since $T$ is compact from $E_{0} \times F_{0}$ into $G_{0}$, it follows, by
Theorem 4.1 that $T \circ (P_{m},Q_{n})$ is compact from $E \times F$ into $G
$.

{\bf Step 2:} For each $m \in \N$, $T \circ (P_{m},I)$ must also be
shown to be compact. We shall show that
\[
\lim_{n \rightarrow \infty} ||T \circ (P_{m},I) - T(P_{m},Q_{n})||_{Bil(E
\times F,G)} = 0.
\]
Since,
\begin{eqnarray*}
T \circ (P_{m},I) - T \circ (P_{m},Q_{n}) & = & T \circ (P_{m}, I - Q_{n}) \\
& = & T \circ (P_{m}, Q^{+}_{n} + Q^{-}_{n}) \\
& = & T \circ(P_{m}, Q^{+}_{n}) + T \circ(P_{m}, Q^{-}_{n}),
\end{eqnarray*}
we need to show that
\begin{equation}
\lim_{n \rightarrow \infty} ||T \circ (P_{m}, Q^{+}_{n})||_{Bil(E \times
F,G)} = \lim_{n \rightarrow \infty} ||T \circ (P_{m}, Q^{-}_{n})||_{Bil(E
\times F,G)} = 0.
\end{equation}
But,
\begin{eqnarray*}
||T \circ (P_{m}, Q^{+}_{n})||_{Bil(E \times F,G)} & \leq & C \, ||T (P_{m}, Q^{+}_{n})||_{0} \, \ov{\rho}\left(\frac{||T (P_{m}, Q^{+}_{n})||_{1}}{||T (P_{m}, Q^{+}_{n})||_{0}}\right) \\
& \leq & C \, \ov{\rho}(||T||_{1}) \, \frac{\ov{\rho}\left(\frac{1}{||T (P_{m}, Q^{+}_{n})||_{0}}\right)}{\frac{1}{||T (P_{m},Q^{+}_{n})||_{0}}}.
\end{eqnarray*}

Hence, it's enough to show that $||T (P_{m},Q^{+}_{n})||_{Bil(E_{0}
\times F_{0},G_{0})} \rightarrow 0$, as $n \rightarrow \infty$.\newline
By contradiction, let us suppose the contrary.\newline Let
$\{(a_{n}, b_{n})\}$ be a bounded sequence in $E_{0} \times F_{0}$.
Since $\{P_{m}\}$ and $\{Q^{+}_{n}\}$ are uniformly bounded in
$E_{0} \times F_{0}$, there is a subsequence $\{n^{\prime}\}$ and a
$\lambda \neq 0$ such that
\[
||T (P_{m} a_{n^{\prime}}, Q^{+}_{n^{\prime}}
b_{n^{\prime}})||_{G_{0}} \rightarrow \lambda \;, \;\;{\rm as} \;\;
n^{\prime}\rightarrow \infty.
\]
By the compactness assumption on $T : E_{0} \times F_{0} \rightarrow
G_{0}$ we may assume, passing the another subsequence if necessary,
that $\{T (P_{m}a_{n^{\prime}} , \; Q^{+}_{n^{\prime}}
b_{n^{\prime}})\}$ converges to some element $b$ in $G_{0}$, so that
$||b||_{G_{0}} = \lambda$. But
\begin{eqnarray*}
||T(P_{m} a_{n^{\prime}}, \; Q^{+}_{n^{\prime}}
b_{n^{\prime}})||_{G_{0} + G_{1}} & \leq & C \, ||P_{m}
a_{n^{\prime}}||_{E_{0} + E_{1}}
||Q^{+}_{n^{\prime}} b_{n^{\prime}}||_{F_{0} + F_{1}} \\
& \leq & C \, ||P_{m} a_{n^{\prime}}||_{E_{0}} ||Q^{+}_{n^{\prime}}
b_{n^{\prime}}||_{F_{1}} \\
& \leq & C \, ||a_{n^{\prime}}||_{E_{0}}
||Q^{+}_{n^{\prime}}||_{L(F_{1},F_{0})} ||b_{n^{\prime}}||_{F_{0}}
\end{eqnarray*}
Since $\lim_{n^{\prime}\rightarrow \infty}
||Q^{+}_{n^{\prime}}||_{L(F_{1},F_{0})} = 0$, it follows that
$T(P_{m} a_{n^{\prime}}, Q^{+}_{n^{\prime}} b_{n^{\prime}})
\rightarrow 0$ in $G_{0} + G_{1}$, as $n \rightarrow \infty$. Consequently $b = 0$, and $\lambda = 0$, which is not the
case.

To prove the second limit in (21) is
also zero, we take a bounded sequence in $E_{1} \times F_{1}$ and
proceed analogously to the first case, but now using the compactness
assumption from $E_{1} \times F_{1} $ into $G_{1}$. Observe that
this assumption can be avoided if we assume that $E_{0}
\hookrightarrow E_{1}$ and $F_{1} \hookrightarrow F_{0}$.

Analogously we prove that $T \circ (I,Q_{n})$ is compact.

{\bf Step 3:} Finally, (ii) shall be proved. We have
\begin{eqnarray}
\lefteqn{T - T \circ(P_{m},I) - T \circ(I,Q_{n}) + T \circ (P_{m},Q_{n}) =}
\\
& = & T \circ (P^{+}_{m} + P^{-}_{m}, Q^{+}_{n} + Q^{-}_{n})  \nonumber \\
& = & T \circ (P^{+}_{m}, Q^{+}_{n}) + T \circ (P^{+}_{m}, Q^{-}_{n}) + T
\circ (P^{-}_{m}, Q^{+}_{n}) + T \circ (P^{-}_{m}, Q^{-}_{n}) .  \nonumber
\end{eqnarray}
To conclude the proof we have to verify that each term on the right
side of (22) converges to zero as $n \rightarrow \infty$. However, in the
proof of (21) the only property used of $P_{m}$ is boundedness. Thus, the
same proof works with $P_{m}$ replaced by $P^{+}_{m}$ or $P^{-}_{m}$. Hence,
$T \circ (P^{+}_{m},Q^{+}_{n})$, $T \circ (P^{-}_{m},Q^{+}_{n})$, $T \circ
(P^{+}_{m},Q^{-}_{n})$ and $T \circ (P^{-}_{m},Q^{-}_{n})$ converge to zero,
as $m,n \rightarrow \infty$.

The proof is complete.

Now, our main goal will be dealt with. We shall state a
bilinear version of Hayakawa's compactness theorem. The idea is to
reduce it to Theorem 5.1.

{\bf Theorem 5.2.} Let $\E = (E_{0}, E_{1}) \; ,
\; \F = (F_{0}, F_{1}) \,$ and $\G = (G_{0}, G_{1})$ be
Banach couples. Let $\, T \in Bil(\E \times \F, \G)$
be given, such that the restrictions $T|_{E_{k} \times F_{k}}$ are
compact from $E_{k} \times F_{k}$ into $G_{k}, \; k = 0, 1$.
Then, given $\rho \in {\cal B}^{+-}$, $T$ is compact from $\E_{\gamma,p} \times \F_{\rho,q}$
into $\G_{\rho,r} \,$, where $\gamma(t) = 1/\ov{\rho}(t^{-1})$ and $\, 1/r = 1/p + 1/q - 1$.

{\bf Proof.} To prove that the bounded bilinear mapping %
\begin{equation}
T : \E_{\gamma,p} \times \F_{\rho,q} \longrightarrow
\G_{\rho,r} \nonumber
\end{equation}
is compact, it is enough to show that
\begin{equation}
\widetilde{T} = T \circ (\sigma, \sigma) : (\ell^{q}_{0}(\Delta_m \E), \ell^{q}_{1}(\Delta_m \E))_{\gamma, q} \times
(\ell^{q}_{0}(\Delta_m \F), \ell^{q}_{1}(\Delta_m \F))_{\rho, q}
\longrightarrow \G_{\rho r} \nonumber
\end{equation}
is compact. But, since the mappings
\[
\ell^{1}_{k} (\Delta_{m} \E) \times \ell^{1}_{k} (\Delta_{m} \F)
\stackrel{(\sigma,\sigma)}{\longrightarrow} E_{k} \times F_{k} \stackrel{T}{\longrightarrow} G_{k}
\;\;\;( k = 0, 1)
\]
are compact, we need to have at hand approximations
sequences which satisfy (AP) so that Theorem 5.1 may be applied. Assuming the existence of such approximations
sequences, it follows by Theorem 4.1 that
\[
\widetilde{T} : (\ell^{1}_{0} (\Delta_{m} \E) , \ell^{1}_{1} (\Delta_{m}
\E))_{\gamma,p} \times (\ell^{1}_{0} (\Delta_{m} \F), \ell^{1}_{1}
(\Delta_{m} \F))_{\rho,q} \longrightarrow (G_{0}, G_{1})_{\rho,r}
\]
is also compact. Finally, taking into account that $(E_{0}, E_{1})_{\rho,q;J} = \ell^{q}_{\gamma}(\Delta_{m}) / {\sigma}^{-1}{(0)}$, it follows that the mapping
$T$ is compact.

It only remains to verify that the Banach couples $%
(\ell^{1}_{0}(\Delta_{m}X), \ell^{1}_{1}(\Delta_{m}X))$, where either $X = \E$
or $X = \F$, satisfy the approximation hypothesis (AP).

For each $n \in \N$, let us consider the cutting
operators $\, P_{n} , P^{+}_{n} \,$ and $\, P^{-}_{n} \,$, defined on $\,
\ell^{1}_{0}(\Delta_{m}) + \ell^{1}_{m}(\Delta_{m}) \,$ by
\begin{eqnarray*}
P_{n} ({u_{m}}) & = & \{\cdots, 0, 0, u_{-n}, u_{-n+1} \cdots, u_{0},\cdots,
u_{n-1}, u_{n}, 0, 0,\cdots \}, \\
P^{+}_{n} ({u_{m}}) & = & \{ \cdots, 0, 0, u_{n+1}, u_{n+2},\cdots \}, \\
P^{-}_{n} ({u_{m}}) & = & \{ \cdots, u_{-n-2}, u_{-n-1}, 0, 0,\cdots \}.
\end{eqnarray*}
We see that $\, I = P_{n} + P^{+}_{n} + P^{-}_{n} \,$ and $P_{n}$, $P_{n}^{+}
$ and $\, P^{-}_{n} \,$ are uniformly bounded in $\, \ell^{1}_{k}
(\Delta_{m}) \; , \; k = 0, 1 \,$. Moreover, $\, P^{+}_{n} : \ell^{1}_{1}
(\Delta_{m}) \rightarrow \ell^{1}_{0} (\Delta_{m}) \,$ and $\, P^{-}_{n} :
\ell^{1}_{0} (\Delta_{m}) \rightarrow \ell^{1}_{1} (\Delta_{m})\,$, and
their norms are bounded by $\, 2^{-n} \,$. Hence the Banach couple $\,
(\ell^{1}_{0}(\Delta_{m}), \ell^{1}_{1} (\Delta_{m})) \,$ verifies the
Approximation Hypothesis (AP).

The proof is complete.

\section{Compactness Theorem of Persson Type}

The first generalization of the Lions-Peetre compactness
theorems was given by A. Persson [14]. The equality between the departure or
arriving spaces is replaced by an approximation hypothesis. The general idea
of Persson's theorem goes back to Krasnoselskii [10], where a compactness
theorem of Riesz-Thorin type is proved.

{\bf Definition 1.} A Banach couple $(E_{0}, E_{1})$ is said to verify
{\em Lions' approximation condition} if a sequence $(P_{n})_{n}$
of linear operators $P_{n}:E_{0} + E_{1} \rightarrow E_{0} + E_{1}$ exists, with $%
P_{n}(E_{k}) \subset E_{0} \cap E_{1}, \;\; k = 0, 1$, and such that $P_{n}
x \rightarrow x$ in $E_{k}$ as $n \rightarrow + \infty$, for each fixed $x
\in E_{k}, \;\; k = 0, 1$.

{\bf Definition 2.} A Banach couple $(E_{0}, E_{1})$ is said to verify
{\em Persson's approximation condition} if, to each compact set $K \subset
E_{0}$, there exists a constant $C > 0$ and a set ${\cal P}$ of linear
operators $P : E_{0} + E_{1} \rightarrow E_{0} + E_{1}$, with $P(E_{k})
\subset E_{0} \cap E_{1}, \;\; k = 0, 1$, such that
\[
||P||_{L(E_{k},E_{k})} \leq C, \;\;\;\;\;\; k = 0, 1.
\]
Furthermore, it is supposed that for each $\varepsilon > 0$ we can find a $P
\in {\cal P}$ so that
\[
||Px - x||_{E_{0}} < \varepsilon,
\]
for all $x \in K$.

{\bf Remark.} Banach--Steinhaus theorem shows that Lions's condition implies Persson's condition.

{\bf Theorem 6.1.} Let $(E_{0}, E_{1})$, $(F_{0}, F_{1})$ and $(G_{0},G_{1})$
be Banach couples such that the pair $(G_{0}, G_{1})$ satisfies Persson's condition.  If $T \in Bil(E_{k}\times F_{k},G_{k})$ is compact from $E_{0} \times F_{0}$ into $%
G_{0}$, given $\rho \in {\cal B}^{+-}$ and $\gamma(t) =1/\ov{\rho}(t^{-1})$, then for spaces $E \in J_\gamma(E_0,E_1) \cap K_\gamma(E_0,E_1)$, $F \in J_\rho(F_0,F_1) \cap K_\rho(F_0,F_1)$ and $G \in J_\rho(G_0,G_1) \cap K_\rho(G_0,G_1)$, $T$ is also compact from $E \times F$ into $G$.

{\bf Proof.} The image $K = T(B_{E_{0}} \times B_{F_{0}})$ in $G_{0}$ of the unit
ball $B_{E_{0}} \times B_{F_{0}}$ of $E_{0} \times F_0$ is relatively compact in $G_{0}$. Hence, choosing $P$
in accordance with Persson's approximation condition, we find
\[
||PT(x,y) - T(x,y)||_{G_{0}} < \varepsilon,
\]
for all $(x,y) \in B_{E_{0}} \times B_{F_{0}}$; that is
\[
||PT - T||_{Bil(E_{0} \times F_{0}, G_{0})} \leq \varepsilon.
\]
From Definition 2, we obtain
\begin{eqnarray*}
||PT - T||_{Bil(E\times F,G)} & \leq & C \, ||PT - T||_{Bil(E_0 \times F_0, G_0)} \; \ov{\rho}\left(\frac{||PT
- T||_{Bil(E_1 \times F_1 , G_1)}}{||PT - T||_{Bil(E_0 \times F_0, G_0)}}\right) \\
& \leq & C \, \ov{\rho}(||T||_{Bil(E_1 \times F_1 , G_1)}) \; \frac{\ov{\rho}\left(\frac{1}{||PT - T||_{Bil(E_0 \times F_0, G_0)}}\right)}{\frac{1}{||PT - T||_{Bil(E_0 \times F_0, G_0)}}}.
\end{eqnarray*}

This means that the bilinear mapping $T : E \times F\rightarrow G$ may be
approximated uniformly by operators of the form $PT$, where $P \in {\cal P}$%
. Hence the theorem will follow if we prove that each mapping $PT : E
\times F \rightarrow G$, with $P \in {\cal P}$, is compact.

According to the closed graph theorem, the mappings $P : G_{k}
\rightarrow G_{0} \cap G_{1}, \;\; k = 0, 1$, are bounded. Since the
composition of a compact and a bounded operator is compact, $PT : E_{0}
\times F_{0} \rightarrow G$ is compact and $PT : E_{1} \times F_{1}
\rightarrow G$ is bounded. Lions--Peetre's Theorem 4.3 shows that $%
PT : E \times F \rightarrow G$ is compact. Proof is thus complete.

\bigskip



\end{document}